\input amstex
\input amsppt.sty
\magnification=\magstep1
\hsize=30truecc
\baselineskip=16truept
\vsize=22.2truecm
\nologo
\pageno=1
\topmatter
\TagsOnRight

\def\up{\bold}
\def\Z{\Bbb Z}

\def\Q{\Bbb Q}

\def\C{\Bbb C}
\def\al{\alpha}
\def\l{\left}
\def\r{\right}
\def\bg{\bigg}
\def\({\bg(}
\def\){\bg)}
\def\[{\bg[}
\def\]{\bg]}
\def\t{\text}
\def\f{\frac}
\def\em{\emptyset}
\def\se {\subseteq}

\def\sm{\setminus}

\def\ls{\leqslant}
\def\gs{\geqslant}

\def\ord{\roman{ord}}
\def\cal{\Cal}
\def\Proof{\noindent{\it Proof}}
\def\Def{\medskip\noindent{\it Definition}\ }
\def\Remark{\medskip\noindent{\it  Remark}}
\def\Ack{\medskip\noindent {\bf Acknowledgment}}
\hbox{Acta Arith. 131(2008), no.\,4, 341--350.}
\bigskip
\title On Covers of Abelian Groups by Cosets\endtitle
\author G\"unter Lettl$^1$ and Zhi-Wei Sun$^2$\endauthor
\leftheadtext{G\"unter Lettl and Zhi-Wei Sun}
\affil $^1$Institut f\"ur Mathematik und wissenschaftliches Rechnen
\\Karl-Franzens-Universit\"at\\ A-8010 Graz, Heinrichstra{\ss}e 36, Austria
\\guenter.lettl\@uni-graz.at
\\{\tt http://www-ang.kfunigraz.ac.at/$\sim$lettl}
\medskip
$^2$Department of Mathematics, Nanjing University\\
Nanjing 210093, People's Republic of China
\\zwsun\@nju.edu.cn
\\{\tt http://math.nju.edu.cn/$\sim$zwsun}
\endaffil

\abstract Let $G$ be any abelian group and $\{a_sG_s\}_{s=1}^k$ be
a finite system of cosets of subgroups $G_1,\ldots,G_k$. We show
that if $\{a_sG_s\}_{s=1}^k$ covers all the elements of $G$ at
least $m$ times with the coset $a_tG_t$ irredundant then
$[G:G_t]\ls 2^{k-m}$ and furthermore $k\gs m+f([G:G_t])$, where
$f(\prod_{i=1}^rp_i^{\al_i})=\sum_{i=1}^r\al_i(p_i-1)$ if
$p_1,\ldots,p_r$ are distinct primes and $\al_1,\ldots,\al_r$ are
nonnegative integers. This extends Mycielski's conjecture in a new
way and implies a conjecture of Gao and Geroldinger. Our new
method involves algebraic number theory and characters of abelian
groups.
\endabstract

\thanks 2000 {\it Mathematics Subject Classification}:\,Primary 20K99;
Secondary 05D99, 05E99, 11B25, 11B75, 11R04, 11S99, 20C15, 20D60.
\newline\indent {\it Key words and phrases}: covers of abelian groups, characters, $p$-adic
valuation.
\newline\indent The preprint was posted as {\tt
arXiv:math.GR/0411144} on Nov. 7, 2004.
\newline\indent The second author is responsible for communications, and supported by
the National Science Fund for Distinguished Young Scholars in
China (grant no. 10425103).
\endthanks
\endtopmatter
\document

\heading{1. Introduction}\endheading

 As in any textbook on group theory, for a subgroup $H$ of a group $G$
 with the index $[G:H]$ finite, $G$ can be partitioned into $k=[G:H]$ left cosets of $H$ in $G$, i.e.,
 all the $k$ left cosets of $H$ form a disjoint cover of $G$.

In 1954 B. H. Neumann [N1, N2] discovered the following basic result on covers of groups.

\proclaim{Theorem 1.1 {\rm (Neumann)}} Let $\{a_sG_s\}_{s=1}^k$ be
a cover of a group $G$ by $($finitely many$)$ left cosets of
subgroups $G_1,\ldots,G_k$. Then $G$ is the union of those
$a_sG_s$ with $[G:G_s]<\infty$. In other words, if
$\{a_sG_s\}_{s\not=t}$ is not a cover of $G$ then we have
$[G:G_t]<\infty$.
\endproclaim

In 1966 J. Mycielski (cf. [MS]) posed an interesting conjecture
on disjoint covers of abelian groups.
Before stating the conjecture we give a definition first.

\Def\ 1.1. The {\it Mycielski function} $f:\Z^+=\{1,2,\ldots\}\to \{0,1,\ldots\}$
is given by
$$f(n)=\sum_{p\in P(n)}\ord_p(n)(p-1),\tag1.1$$
where $P(n)$ denotes the set of prime divisors of $n$ and $\ord_p(n)$
represents the largest nonnegative integer $\al$ such that $p^{\al}\mid n$.
\medskip

\Remark\ 1.1. Since $p\ls 2^{p-1}$ for any prime $p$, (1.1)
implies that $n\ls 2^{f(n)}$ (i.e., $f(n)\gs\log_2 n$).

\medskip
\proclaim{Mycielski's Conjecture} Let $G$ be an abelian group, and $\{a_sG_s\}_{s=1}^k$
be a disjoint cover of $G$ by left cosets of subgroups.
Then $k\gs 1+f([G:G_t])$ for each $t=1,\ldots,k$.
\endproclaim

When $G$ is the additive group $\Z$ of integers, Mycielski's conjecture says that
for any disjoint cover $\{a_s(n_s)\}_{s=1}^k$ of $\Z$ by residue classes
(where $a_s\in\Z$, $n_s\in\Z^+$ and $a_s(n_s)=a_s+n_s\Z$) we have $k\gs1+f(n_t)$
for every $t=1,\ldots,k$. This was first confirmed by \v S. Zn\'am [Z66].
For problems and results on covers of $\Z$, the reader is referred to [G04], [PS], [S03]
and [S05].

\Def\ 1.2. For a subnormal subgroup $H$ of a group $G$ with finite index, we define
$$d(G,H)=\sum_{i=1}^n([H_{i}:H_{i-1}]-1),\tag1.2$$
where $H_0=H\subset H_1\subset \cdots\subset H_n=G$
is any composition series from $H$ to $G$.
\medskip

 By [S90, Theorem 6] and [S01, Theorem 3.1], for any subnormal subgroup $H$ of a group $G$ with $[G:H]<\infty$,
 we have $d(G,H)\gs f([G:H])$, and equality holds if and only if $G/H_G$ is solvable,
 where $H_G=\bigcap_{g\in G}gHg^{-1}$ is the core of $H$ in $G$ (i.e., the largest normal subgroup
 of $G$ contained in $H$).

The following result is stronger than Mycielski's conjecture.

\proclaim{Theorem 1.2 {\rm (I. Korec, Z. W. Sun)}} Let $a_1G_1,\ldots,a_kG_k$ be
left cosets of subnormal subgroups $G_1,\ldots,G_k$ of a group $G$.
If ${\cal A}=\{a_sG_s\}_{s=1}^k$ forms an exact $m$-cover of $G$, i.e.,
${\cal A}$ covers each element of $G$ exactly $m$ times, then
$[G:\bigcap_{s=1}^kG_s]<\infty$ and
$$k\gs m+d\(G,\bigcap_{s=1}^kG_s\)\gs m+f\(\[G:\bigcap_{s=1}^kG_s\]\),$$
where the lower bound $m+d(G,\bigcap_{s=1}^kG_s)$ is best possible.
\endproclaim

In the case $m=1$ and $G=\Z$, Theorem 1.2 was first conjectured by Zn\'am [Z69].
When $m=1$ and $G_1,\ldots,G_k$ are normal in $G$, Theorem 1.2 was obtained by Korec [K74]
in 1974. In 1990 Sun [S90] deduced Theorem 1.2 in the case $m=1$ by a method different
from that of Korec. The current version of Theorem 1.2 was established by Sun [S01] in 2001,
the proof of which depends heavily on the condition that ${\cal A}$ covers all the elements
of $G$ the same number of times. Under the conditions of Theorem 1.2, Sun [S04] also showed that
the indices $[G:G_s]\ (1\ls s\ls k)$ cannot be distinct providing $k>1$.

Call a coset in an abelian group not containing the identity element
a {\it proper coset}. In 2003 W. D. Gao and A. Geroldinger [GG] proved the following
conjecture for any elementary abelian $p$-group $G$ (they did not explicitly state this conjecture in [GG]).

\proclaim{Gao-Geroldinger Conjecture} Let $G$ be a finite abelian group with identity $e$.
If $G\sm\{e\}$ is a union of $k$ proper cosets $a_1G_1,\ldots,a_kG_k$ then
we have $k\gs f(|G|)$.
\endproclaim

With the notations of the Gao-Geroldinger conjecture, if we set $a_0=e$ and $G_0=\{e\}$
then $\{a_sG_s\}_{s=0}^k$ forms
a cover of $G$ with $a_0G_0\cap a_sG_s=\em$ for all $s=1,\ldots,k$. Thus, by the result of [Z69],
the Gao-Geroldinger conjecture holds when $G$ is cyclic.

 In this paper we aim to generalize Mycielski's conjecture in a new direction
 and prove an extended version of the Gao-Geroldinger conjecture.

\medskip
\Def\ 1.3. Let $G$ be a group and let ${\cal A}=\{a_sG_s\}_{s=1}^k$ be
a finite system of left cosets of subgroups $G_1,\ldots,G_k$. The {\it covering function}
of ${\cal A}$ is given by
$$w_{{\cal A}}(x)=|\{1\ls s\ls k:\,x\in a_sG_s\}|\ \ \ (x\in G).\tag1.3$$
Let $m$ be a positive integer. We call ${\cal A}$ an {\it $m$-cover} of $G$ if $w_{{\cal A}}(x)\gs m$
for all $x\in G$. If ${\cal A}$ forms an $m$-cover of $G$ but none of its proper subsystems
does, then ${\cal A}$ is said to be a {\it minimal $m$-cover} of $G$.
\medskip

 Now we state our main result in this paper,
 which (in the special case $m=1$)
 implies the Gao-Geroldinger conjecture for arbitrary finite abelian groups.

\proclaim{Theorem 1.3} Let ${\cal A}=\{a_sG_s\}_{s=1}^k$ be an $m$-cover of an abelian group $G$
by left cosets. Then, for any $a\in G$ with $w_{{\cal A}}(a)=m$,
we have
$$N_a=\[G:\bigcap\Sb 1\ls s\ls k\\a\in a_sG_s\endSb G_s\]\ls 2^{k-m}
\ \ \t{and furthermore}\ \ k\gs m+f(N_a).\tag1.4$$
 In particular, if
$\{a_sG_s\}_{s\not=t}$ fails to be an $m$-cover of $G$, then we have the inequalities
$$[G:G_t]\ls 2^{k-m}\quad\t{and}\quad k\gs m+f([G:G_t]),\tag1.5$$
the bounds of which are best possible.
\endproclaim

\Remark\ 1.2. When $G=\Z$, Theorem 1.3 was proved  by Zn\'am [Z75]
in the case $m=1$, and we can say something stronger in Section 2.
Also, in the second inequality of (1.4), $N_a$ cannot be replaced
by $[G:\bigcap_{s=1}^kG_s]$ as illustrated by the following
example.

\medskip
{\it Example} 1.1. Let $G$ be the abelian group $C_p\times C_p$ where $p$ is a prime
and $C_p$ is the cyclic group of order $p$. Then any element $a\not=e$ of $G$ has order $p$.
 Let $G_1,\ldots,G_k$ be all the distinct
subgroups of $G$ with order $p$.
If $1\ls i<j\ls k$, then $G_i\cap G_j=\{e\}$. Thus
$\{G_s\}_{s=1}^k $ forms a minimal $1$-cover of $G$ with
$\bigcap_{s=1}^kG_s=\{e\}$.
Since $1+k(p-1)=|\bigcup_{s=1}^kG_s|=|G|=p^2$, we have
  $$k=p+1 \gs 1+f([G:G_s])=1+f(p)=p.$$
However,
 $$k=p+1 \ls 2p-1=1+f([G:\{e\}])=1+d\(G,\bigcap^k_{s=1}G_s\),$$
 and the last inequality becomes strict when $p>2$.
\medskip

Example 1.1 also shows that we don't have an analogy of [S01,
Theorem 2.1] for minimal $m$-covers of the abelian group $C_p\times
C_p$ (where $p$ is a prime), thus we cannot prove our Theorem
1.3 by the method in [S01]. To obtain Theorem 1.3 we employ some
tools from algebraic number theory as well as characters of
abelian groups.

\proclaim{Corollary 1.1} Let ${\cal A}=\{a_sG_s\}_{s=1}^k$ be an $m$-cover
of a group $G$ by left cosets. Provided that $a\in G$ and $w_{{\cal A}}(a)=m$,
for any abelian subgroup $K$ of $G$ we have
$$\aligned k-m\gs&|\{1\ls s\ls k:\ a\not\in a_sG_s\ \t{and}\ K\not\se G_s\}|
\\\gs&f\(\[K:K\cap\bigcap^k\Sb s=1\\a\in a_sG_s\endSb G_s\]\).
\endaligned\tag1.6$$
In particular, if $\{a_sG_s\}_{s\not=t}$ fails to be an $m$-cover
of $G$, then for any abelian subgroup $K$ of $G$ not contained in $G_t$ we
have
$$|\{1\ls s\ls k:\, K\not\se G_s\}|\gs1+f([K:G_t\cap K]).\tag1.7$$
\endproclaim
\Proof. Let $J=\{1\ls s\ls k:\, a_sG_s\cap aK\not=\em\}$.
For each $s\in J$, $a^{-1}a_sG_s\cap K$ is a coset of $G_s\cap K$ in $K$.
Observe that $\{a^{-1}a_sG_s\cap K\}_{s\in J}$ is an $m$-cover of $K$
with $|\{s\in J:\,e\in a^{-1}a_sG_s\cap K\}|=|I_a|=m$
where $I_a=\{1\ls s\ls k:\,a\in a_sG_s\}$.
Applying Theorem 1.3 to the abelian group $K$ we get the inequality
$|J|-m\gs f([K:\bigcap_{s\in I_a}G_s\cap K])$.
If $s\in J$ and $K\se G_s$, then $a^{-1}a_sG_s\cap K=K$ and hence $s\in I_a$.
Thus
$$\aligned |J|-m=&|\{s\in J:\,e\not\in a^{-1}a_sG_s\cap K\}|
\\\ls&|\{1\ls s\ls k:\, a\not\in a_sG_s\ \t{and}\ K\not\se G_s\}|
\ls k-m\endaligned$$
and hence (1.6) follows.

Now suppose that $\{a_sG_s\}_{s\not=t}$ is not an
$m$-cover of $G$ and $K$ is an abelian subgroup of $G$ with $K\not\se G_t$.
Then $w_{{\cal A}}(x)=m$ for some $x\in a_tG_t$. In light of the above,
$$\aligned &|\{1\ls s\ls k:\, s\not=t\ \t{and}\ K\not\se G_s\}|
\\\gs&|\{1\ls s\ls k:\, x\not\in a_sG_s\ \t{and}\ K\not\se G_s\}|\gs f([K:K\cap G_t]).
\endaligned$$
This proves (1.7) and we are done. \qed

\proclaim{Corollary 1.2} Let $R$ be any ring. Let $a_1,\ldots,a_k$ be elements of $R$ and
$I_1,\ldots,I_k$ ideals of $R$. If $\{a_s+I_s\}_{s=1}^k$
is an $m$-cover of $R$ with the coset $a_t+I_t$ irredundant,
then for the quotient ring $R/I_t$ we have $|R/I_t|\ls 2^{k-m}$ and furthermore $k\gs m+f(|R/I_t|)$.
\endproclaim
\Proof. Since $R$ is an additive abelian group, this follows from Theorem 1.3 immediately. \qed

In the next section we will present a new approach to Mycielski's
problem on covers of $\Z$. In Section 3 we are going to work with
covers of abelian groups and extend some ideas in Section 2, this
will lead to our proof of Theorem 1.3.

\heading{2. A new approach to Mycielski's problem}\endheading

Let $\overline{\Q}$ denote the algebraic closure of the rational field $\Q$ and
$\overline{\Z}$ the ring of all algebraic integers in $\overline{\Q}$.

\proclaim{Lemma 2.1} For $s=1,\ldots,k$ let $\zeta_s\in\overline{\Z}$ be a
root of unity with order $n_s>1$. Then $n\in\Z^+$ divides $\prod_{s=1}^k(1-\zeta_s)$ in
$\overline{\Z}$, if and only if we have
$$\sum^k\Sb
s=1\\P(n_s)=\{p\}\endSb\f1{\varphi(n_s)}\gs\ord_p(n)\quad \t{for any prime}\ p,\tag2.1$$
where $\varphi$ is the well-known Euler function.
\endproclaim
\Proof. For each prime $p$, let $\up{v}_p:\overline{\Q}\to\Q$
denote any extension of the $p$-adic valuation $\ord_p(\cdot)$ to
$\overline{\Q}$, normed by $\up{v}_p(p)=1$. It is well known (cf.
[W, Chap. 2]) that
$$\up{v}_p(1-\zeta_s)=\cases1/\varphi(n_s)&\t{if}\ n_s\ \t{is a power of}\ p,
\\0&\t{otherwise.}\endcases$$
Now $n$ divides $\prod_{s=1}^k(1-\zeta_s)$ in $\overline{\Z}$, if
and only if for each valuation $\up{v}:\overline{\Q}\to\Q$ one has
$\up{v}(n)\ls\sum_{s=1}^k\up{v}(1-\zeta_s)$. Since any valuation
$\up{v}$ of $\overline{\Q}$ is (equivalent to) an extension of
$\ord_p(\cdot)$ for some prime $p$, we immediately obtain the
desired result. \qed

\proclaim{Corollary 2.1} Let $n>1$ be an integer. Then $f(n)$ is
the smallest positive integer $k$ such that there are roots of unity
$\zeta_1,\ldots,\zeta_k$ different from $1$ for which
$\prod_{s=1}^k(1-\zeta_s)\in n\overline{\Z}$.
Furthermore, this holds with $k=f(n)$ if and only if for any prime divisor $p$ of $n$ there
are exactly $\ord_p(n)(p-1)$ of $\zeta_1,\ldots,\zeta_k$ having order $p$.
\endproclaim

\Proof. For $s=1,\ldots,k$ let $\zeta_s$ be a root of unity with order $n_s>1$.
By Lemma 2.1, $n$ divides $\prod_{s=1}^k(1-\zeta_s)$ in $\overline{\Z}$ if and only if (2.1) holds.
Clearly
$$\sum^k\Sb s=1\\P(n_s)=\{p\}\endSb\f1{\varphi(n_s)}\ls \f{|\{1\ls s\ls k:\, P(n_s)=\{p\}\}|}{p-1}
\ \ \ \ \t{for every prime}\ p.$$
If (2.1) is valid, then
$$k\gs\sum_{p\in P(n)}|\{1\ls s\ls k:\, P(n_s)=\{p\}\}|
\gs\sum_{p\in P(n)}\ord_p(n)(p-1)=f(n).$$

 Now assume that $k=f(n)$. When (2.1) is valid, equality holds in the last three inequalities and hence
$$|\{1\ls s\ls k:\,n_s=p\}|=|\{1\ls s\ls k:\,P(n_s)=\{p\}\}|=\ord_p(n)(p-1)$$
for any prime $p$. Conversely, (2.1) holds if $|\{1\ls s\ls k:\,n_s=p\}|=\ord_p(n)(p-1)$
for all $p\in P(n)$.

 Combining the above we have completed the proof. \qed

\proclaim{Lemma 2.2} Suppose that $A=\{a_s(n_s)\}_{s=1}^k$ is an
$m$-cover of $\Z$ by residue classes and $a\in\Z$
is covered by $A$ exactly $m$ times.
Let $N_a$ be
the least common multiple of those $n_s$ with $a\in a_s(n_s)$,
and let $m_s\in\Z$ for $s\in J$ where $J=\{1\ls s\ls k:\, a\not\in a_s(n_s)\}$.
Then, for any $0\ls\al<1$ we have
$$C_0(\al)=C_1(\al)=\cdots=C_{N_a-1}(\al),\tag2.2$$
where
$$C_r(\al)=\sum\Sb I\se J\\\{\sum_{s\in I}m_s/n_s\}=(\al+r)/N_a\endSb
(-1)^{|I|}e^{2\pi i\sum_{s\in I}(a_s-a)m_s/n_s}\tag2.3$$
for every $r=0,1,\ldots,N_a-1$,
and we use $\{\theta\}$ to denote the fractional part of a real number $\theta$.
\endproclaim
\Proof. This follows from [S99, Lemma 2]. \qed

\proclaim{Theorem 2.1} Let $A=\{a_s(n_s)\}_{s=1}^k$ be an
$m$-cover of $\Z$, and suppose that $a$ is an integer with $w_A(a)=m$. Then
$k\gs m+f(N_a)$ where $N_a$ is the least common multiple of
those $n_s$ with $a\in a_s(n_s)$.
Furthermore, for any prime $p$ we have
$$|I(p)|\gs\sum_{s\in I(p)}\f1{p^{\ord_p(n_s)-\ord_p(a_s-a)-1}}\gs \ord_p(N_a)(p-1),\tag2.4$$
where
$$I(p)=\bg\{1\ls s\ls k:\ \f{n_s}{p^{\ord_p(n_s)}}\ \bg|\  a_s-a\ \t{but}\ n_s\nmid a_s-a\bg\}.\tag2.5$$
\endproclaim
\Proof. Let $J=\{1\ls s\ls k:\, a\not\in a_s(n_s)\}$. For each $s\in J$,
let $m_s$ be an integer
not divisible by $n_s/(n_s,a_s-a)>1$, then $\zeta_s=e^{2\pi i(a_s-a)m_s/n_s}$ is
a primitive $d_s$th root of unity where $d_s=n_s/(n_s,(a_s-a)m_s)>1$.

 Set
$$S=\bg\{\bg\{N_a\sum_{s\in I}\f{m_s}{n_s}\bg\}:\ I\se J\bg\}.$$
Then
$$\aligned\prod_{s\in J}\l(1-\zeta_s\r)=&\sum_{I\se J}(-1)^{|I|}
e^{2\pi i\sum_{s\in I}(a_s-a)m_s/n_s}
\\=&\sum_{\al\in S}\sum\Sb I\se J\\\{N_a\sum_{s\in I}m_s/n_s\}=\al\endSb
 (-1)^{|I|}e^{2\pi i\sum_{s\in I}(a_s-a)m_s/n_s}
\\=&\sum_{\al\in S}\sum_{r=0}^{N_a-1}C_r(\al)=N_a\sum_{\al\in S}C_0(\al),
\endaligned$$
where $C_r(\al)\ (0\ls r<N_a)$ are given by (2.3). So
$N_a$ divides $\prod_{s\in J}(1-\zeta_s)$ in the ring $\overline{\Z}$.
By Corollary 2.1, we have $k-m=|J|\gs f(N_a)$.
In view of Lemma 2.1,
$$\sum\Sb s\in J\\P(d_s)=\{p\}\endSb\f1{\varphi(d_s)}\gs\ord_p(N_a)
\ \ \t{for each prime}\ p.$$

Now we simply let $m_s=1$ for all $s\in J$.
By the above, for any prime $p$ we have
$$\sum_{s\in I(p)}\f1{\varphi(n_s/(n_s,a_s-a))}\gs \ord_p(N_a)$$
which is equivalent to (2.4).
This concludes the proof. \qed

\heading{3. Working with abelian groups}\endheading

We first recall some well-known facts from the theory of characters of finite
abelian groups (see, e.g. [W, pp.\,22-23]).

For a finite abelian group $G$, let $\widehat G$ denote the group of
all complex-valued characters of $G$. One has $\widehat G\cong G$.
For any subgroup $H$ of $G$ let $H^\bot$ denote the group of those characters
$\chi\in\widehat G$ with $\ker(\chi)=\{x\in G:\,\chi(x)=1\}$
containing $H$. Then there is a canonical isomorphism $H^\bot\cong\widehat{G/H}$ by
putting $\chi(aH)=\chi(a)$ for any $a\in G$ and any $\chi\in H^\bot$.
Furthermore, for each $a\in G\sm H$ there exists some $\chi\in H^\bot$ with
$\chi(a)\not=1$.

\medskip
\noindent{\tt Proof of Theorem 1.3}. Choose a minimal
$I_*\se\{1,\ldots,k\}$ such that the system $\{a_sG_s\}_{s\in I_*}$ forms an
$m$-cover of $G$. As $I_a=\{1\ls s\ls k:\, a\in a_sG_s\}$ has
cardinality $m$, $I_a$ is contained in $I_*$. So we can simply assume that
${\cal A}$ is a minimal $m$-cover of $G$ (i.e.,
$I_*=\{1,\ldots,k\}$). By [S90, Corollary 1],
$H=\bigcap_{s=1}^kG_s$ is of finite index in $G$. Instead of the
minimal $m$-cover ${\cal A}=\{a_sG_s\}_{s=1}^k$ of $G$, we may
consider the minimal $m$-cover $\bar{\cal A}=\{\bar a_s\bar
G_s\}_{s=1}^k$ of the finite abelian group $\bar G=G/H$, where
$\bar a_s=a_sH$ and $\bar G_s=G_s/H$ (hence $[\bar G:\bar G_s]=[G:G_s]$).
Therefore, without any loss of generality, we can assume that $G$
is finite.

Put $H_a=\bigcap_{s\in I_a}G_s$; then $|H_a^\bot|=[G:H_a]=N_a$.

Note that $J=\{1\ls j\ls k:\, a\not\in a_jG_j\}$ has cardinality
$k-m$. For each $j\in J$ we may choose a $\chi_j\in G_j^\bot$
with $\zeta_j:=\chi_j(a^{-1}a_j)\not=1$.
For any $x\in G\sm H_a$ we have $ax\not\in\bigcap_{s\in I_a}aG_s=\bigcap_{s\in I_a}a_sG_s$.
Since ${\cal A}$ is an $m$-cover of $G$, there exists some $j\in J$ with $ax\in a_jG_j$, and therefore
$\chi_j(x)=\zeta_j$
by the choice of $\chi_j$ and the definition of $\zeta_j$.

For $x\in G$ we define
$$\Psi(x)=\prod_{j\in J}\l(\chi_j(x)-\zeta_j\r).$$
If $\chi\in H_a^\bot$ and $\chi(x)\not=1$, then $x\not\in H_a$ and hence
$\Psi(x)=0$ by the above. Thus $\Psi\chi=\Psi$ for all $\chi\in H_a^\bot$.

Observe that
$$\Psi(x)=\sum_{I\se J}\(\prod_{j\in I}\chi_j(x)\)\prod_{j\in J\sm I}\l(-\zeta_j\r)
=\sum_{\psi\in \widehat G}c(\psi)\psi(x),$$
where
$$c(\psi)=\sum\Sb I\se J\\\prod_{j\in I}\chi_j=\psi\endSb\prod_{j\in J\sm I}
\l(-\zeta_j\r)\in\overline{\Z}.$$
Let $\C$ be the complex field. As the set $\widehat G$ is a basis of the $\C$-vector space
$$\C^G=\{g:\ g\ \t{is a function from}\ G\ \t{to}\ \C\}$$
(cf. [J, p.\,291]), for any $\chi\in H_a^\bot$ we have
$c(\psi\chi)=c(\psi)$ for all $\psi\in\widehat G$ because
$\Psi\chi^{-1}=\Psi$.

Clearly
 $$\prod_{j\in J}\l(1-\zeta_j\r)=\Psi(e)=\sum_{\psi\in\widehat G}c(\psi)\psi(e)
 =\sum_{\psi\in\widehat G}c(\psi).$$
 Let $\psi_1H_a^\bot\cup\cdots\cup\psi_lH_a^\bot$ be a coset decomposition of $\widehat G$ where
 $l=[\widehat G:H_a^\bot]$. Then
 $$\sum_{\psi\in\widehat G}c(\psi)=\sum_{r=1}^l\sum_{\chi\in H_a^\bot}c(\psi_r\chi)
 =\sum_{r=1}^l|H_a^\bot|c(\psi_r)=N_a\sum_{r=1}^lc(\psi_r).$$
 (That $c(\psi_r\chi)=c(\psi_r)$ for all $\chi\in H_a^\bot$ is an analogy of Lemma 2.2.)
 Therefore $N_a$ divides $\prod_{j\in J}(1-\zeta_j)$ in $\overline{\Z}$, and
 Corollary 2.1 gives $k-m=|J|\gs f(N_a)$, and consequently $N_a\ls 2^{k-m}$ by Remark 1.1.

 If $\{a_sG_s\}_{s\not=t}$ is not an $m$-cover of $G$, then
 for some $x\in a_tG_t$ we have $w_{{\cal A}}(x)=m$, hence
 $k-m\gs f(N_x)\gs f([G:G_t])$ and $[G:G_t]\ls N_x\ls 2^{k-m}$ by the above.

 By [S01, Example 1.2], for any subgroup $H$ of $G$ (with $[G:H]<\infty$) and an arbitrary element $x$ of $G$,
 the coset $xH$ and $m-1+d(G,H)=m-1+f([G:H])$ other cosets of subgroups containing $H$
 form an (exact) $m$-cover of
 $G$ with $xH$ irredundant. Also, $m-1$ copies of $0(1)$,
together with the following $k-m+1$ residue classes
$$1(2),\ 2(2^2),\ \ldots,\ 2^{k-m-1}(2^{k-m}),\ 0(2^{k-m}),$$
clearly form an (exact) $m$-cover of $\Z$ with the residue class $0(2^{k-m})$ irredundant.
So the inequalities in (1.5) are really best possible and we are done. \qed

\Ack. The authors met each other during the second author's visit
to Graz University in June 2004, so the second author wishes to
thank Prof. A. Geroldinger for the invitation and hospitality.


\widestnumber\key{K74}
\Refs

\ref\key GG\by W. D. Gao and A. Geroldinger, {\it
Zero-sum problems and coverings by proper cosets},
European J. Combin. {\bf 24}(2003), 531--549\endref

\ref\key G04\by R. K. Guy, {\it Unsolved Problems in Number Theory}, 3rd ed.,
Springer, 2004, Sections F13 and F14\endref

\ref\key J\by N. Jacobson, {\it Basic Algebra II}, 2nd ed.,
Freeman \& Co., 1985\endref

\ref\key K74\by I. Korec, {\it On a generalization of Mycielski's
and Zn\'am's conjectures about coset
decomposition of Abelian groups},
Fund. Math. {\bf 85}(1974), 41--48\endref

\ref\key MS\by J. Mycielski and W. Sierpi\'nski,
{\it Sur une propri\'et\'e des ensembles lin\'eaires},
Fund. Math. {\bf 58}(1966), 143--147\endref

\ref\key N1\by B. H. Neumann, {\it Groups covered by permutable subsets},
J. London Math. Soc. {\bf 29}(1954), 236--248\endref

\ref\key N2\by B. H. Neumann, {\it Groups covered by finitely many cosets},
Publ. Math. Debrecen {\bf 3}(1954), 227--242\endref

\ref\key PS\by \v S. Porubsk\'y and J. Sch\"onheim, {\it Covering
systems of Paul Erd\"os: past, present and future}, in: Paul
Erd\"os and his Mathematics. I (edited by G. Hal\'asz, L.
Lov\'asz, M. Simonvits, V. T. S\'os), Bolyai Soc. Math. Studies
11, Budapest, 2002, pp. 581--627\endref

\ref\key S90\by Z. W. Sun, {\it  Finite coverings of groups},
Fund. Math. {\bf 134}(1990), 37--53\endref

\ref\key S99\by Z. W. Sun, {\it  On covering multiplicity},
Proc. Amer. Math. Soc. {\bf 127}(1999), 1293--1300\endref

\ref\key S01\by  Z. W. Sun, {\it  Exact $m$-covers of groups by cosets},
European J. Combin. {\bf 22}(2001), 415--429\endref

\ref\key S03\by  Z. W. Sun, {\it Unification of zero-sum problems,
subset sums and covers of $\Z$},
Electron. Res. Announc. Amer. Math. Soc. {\bf 9}(2003), 51--60\endref

\ref\key S04\by  Z. W. Sun, {\it On the Herzog-Sch\"onheim conjecture for uniform covers of
groups}, J. Algebra {\bf 273}(2004), 153--175\endref

\ref\key S05\by  Z. W. Sun, {\it On the range of a covering function},
J. Number Theory {\bf 111}(2005), 190--196\endref

\ref\key W\by L. C. Washington, {\it Introduction to Cyclotomic Fields},
Springer, New York, 1982\endref

\ref\key Z66\by \v S. Zn\'am, {\it  On Mycielski's problem on systems
of arithmetical progressions},
Colloq. Math. {\bf 15}(1966), 201--204\endref

\ref\key Z69\by \v S. Zn\'am, {\it  A remark to a problem of J. Mycielski
on arithmetic sequences},
Colloq. Math. {\bf 20}(1969), 69--70\endref

\ref\key Z75\by \v S. Zn\'am, {\it  On properties of systems
of arithmetic sequences},
Acta Arith. {\bf 26}(1975), 279--283\endref

\endRefs
\enddocument